\documentclass{amsart}
\usepackage[T1]{fontenc}
\usepackage[latin1]{inputenc}
\usepackage[cm]{aeguill}
\usepackage[francais]{babel}
\usepackage{amsmath,amssymb,enumerate}

\newtheorem{prop}{Proposition}[subsection]
\newtheorem{coro}[prop]{Corollaire}
\newtheorem{lem}[prop]{Lemme}

\newtheorem{conj}{Conjecture}

\newtheorem{ques}[conj]{Question}
\newtheorem{theoi}{Theor\`eme}
\newtheorem{coroi}[theoi]{Corollaire}

\theoremstyle{definition}

\newcommand{\PP}{\mathbb {P}^{n+1}}
\newcommand{\R}{\mathbb R}
\newcommand{\Q}{\mathbb Q}
\newcommand{\C}{\mathbb C}
\newcommand{\N}{\mathbb N}
\newcommand{\V}{\mathbb V}

\newcommand{\Z}{\mathbb Z}

\newcommand{\OSp}{{\mathcal O}_{\Sigma,P}}
\newcommand{\OS}{{\mathcal O}_{\Sigma}}
\newcommand{\OP}{{\mathcal O}_{\mathbb {P}^{n+1}}}

\newcommand{\hdot}{\bullet}

\newcommand{\OX}{{\mathcal O}_X}
\newcommand{\OmX}{\omega_{\mathbb{P}^{n+1}}}
\newcommand{\OXp}{{\mathcal O}_{X,P}}
\newcommand{\G}{\mathcal{G}}
\newcommand{\wt}{\widetilde}

\newenvironment{prv}{{\bf Preuve:}}{$\Box$}

\begin{document}
\title[Application des p\'eriodes des hypersurfaces \`a singularit\'es simples]{Sur l'application des p\'eriodes d'une Variation de Structure de Hodge attach\'ee 
aux familles d'hypersurfaces \`a singularit\'es simples. }
\author{ Philippe Eyssidieux,  Damien M\'egy}
\date{16 Mai 2013}
\begin{abstract}
Soit $n\in \N^*$ un entier positif pair et $d$ un entier positif. Pour toute famille compl\`ete 
$Z$ d'hypersurfaces de $\mathbb{P}^{n+1}$ de degr\'e $d$ \`a singularit\'es isol\'ees de type A-D-E, nous construisons
d'apr\`es une id\'ee de Carlson et Toledo reprise dans \cite{Sim4,Megy}
un champ de Deligne-Mumford $\bar Z$ d'espace de modules $Z$ auquel la repr\'esentation de monodromie de la famille
se prolonge. Nous \'etudions l'application de p\'eriodes associ\'ee et montrons un th\'eor\`eme de Torelli infinit\'esimal 
le long des strates isosinguli\`eres de $Z$ sous des hypoth\`eses de transversalit\'e. Enfin, nous appliquons
ce r\'esultat \`a l'\'etude du rev\^etement universel de $\bar Z$. 
\end{abstract}

\maketitle
 
Dans \cite{Sim4} est d\'ecrite une classe de surfaces projectives alg\'ebriques $S$
munies de
$\Q$-Variations de Structure de Hodge $\V_S$ qui sont int\'eressantes
du point de vue de la th\'eorie de Hodge non-ab\'elienne: $(S,\V_S)$ ne 
peut pas s'exprimer par tir\'e en arri\`ere \`a partir
 de syst\`emes locaux sur des courbes, vari\'et\'es ab\'eliennes
ou espaces localement sym\'etriques hermitiens.
Ce sont des exemples particuli\`erement
int\'eressants pour l'uniformisation en plusieurs variables
complexes (voir \cite{pasy} pour un survey r\'ecent) et l'un de nous a g\'en\'eralis\'e cette construction jusqu'en dimension $6$ et
a entam\'e l'\'etude cohomologique de ces exemples \cite{Megy}. La conjecture de Toledo
stipulant que $H^2(\pi_1(S),\Q)\not =0$ 
  n'est d\'ecid\'ee dans cette classe d'exemples
 que dans certains cas \cite{Megy}. 
La motivation initiale de ce travail
est d'\'etudier  pour cette classe d'exemples l'autre probl\`eme ouvert g\'en\'eral de l'uniformisation
en plusieurs variables complexes, c'est \`a dire la 
conjecture de Shafarevich pr\'edisant que le rev\^etement universel d'une vari\'et\'e
projective alg\'ebrique complexe est holomorphiquement convexe (cf.   
\cite{pasy} pour la d\'efinition de la convexit\'e holomorphe et une discussion du probl\`eme). 

D\'ecrivons la construction de \cite{Sim4,Megy} qui reprend  une id\'ee de Carlson et Toledo. 
Dans ce qui suit $X$ d\'esigne une vari\'et\'e projective complexe  connexe
 de dimension $n+1$, $n\ge 1$, $L$ un faisceau inversible tel que $|L|$ n'a pas de point base. 
Soient $\mathcal{X} \subset |L| \times X$ l'hypersurface universelle
et $p_1: \mathcal{X} \to |L|$ la projection sur le premier facteur.
Si $f\in H^0(X,L)-\{ 0 \}$, on note $X_f=\{ x\in X \ | \  f(x)=0 \}$ et $[f]\in |L|$ le point correspondant. 
D\'efinissons l'ouvert de Zariski  $U(0):=U(X,L)(0)$ comme le lieu
des $[f]\in |L|$ tels que $X_f$ est une hypersurface lisse et notons $D:=|L|-U(0)$ le lieu 
discriminant. On pose $\mathcal{X}(0)=p_1^{-1} (U(0))$. 
Sur $U(0)$, on construit le syst\`eme local de monodromie
$R^n p_{1*} \Q_{\mathcal{X}(0)}$
\footnote{Pour tout espace $T$ et tout groupe $A$,  $A_T$ d\'esigne le faisceau des fonctions localement
constantes sur $T$ \`a valeurs dans $A$. Plus g\'en\'eralement si $\mathbb{W}_T$ est un syst\`eme local sur $T$ et $\phi:T'\to T$
une application continue,  on note $\mathbb{W}_{T'}=\phi^*\mathbb{W}_{T'}$. De m\^eme, pour $X\to T$ une application continue,
on note $X_{T'}=X\times_T T'$. }.
On d\'esigne par
$\V_{U(0)}$ le conoyau du
morphisme $H^n(X, \Q) \otimes \Q_{U(0)} \to R^n p_{1*} \Q_{\mathcal{X}(0)}$.
On a, par le th\'eor\`eme de semi-simplicit\'e,
 un isomorphisme $R^n p_{1*} \Q_{\mathcal{X}(0)}
\simeq  \V_{U(0)} \oplus H^n(X, \Q) \otimes \Q_{U(0)}$. Fixons une fois pour 
toute un \'el\'ement g\'en\'eral $f_{gen}\in H^0(X,L)$ et prenons $\eta:=[f_{gen}]\in U(0)$ comme  point base.

Supposons  d\'esormais que $n$ est pair. Le syst\`eme local $\V_{U(0)}$ est le produit tensoriel par $\Q$ d'un syst\`eme 
sous jacent \`a une $\Z$-Variation de Structures de Hodge de poids $n$ dont la polarisation
est orthogonale et la fibre en $\eta$ est $H_{ev}^n(X_{\eta}, \Q) := \ker(H^n(X_{\eta}, \Q)\to H^n(X, \Q))$. 
On note par $$\rho:\pi_1(U(0), \eta) \to O (H_{ev}^n(X_{\eta}, \Q), \int_{X_{\eta}} - \cup -)$$
la repr\'esentation d'holonomie du syst\`eme local $\V_{U(0)}$.

Introduisons $U:=U(X,L) \subset |L|$ l'ouvert de Zariski form\'e des 
hypersurfaces n'ayant que de singularit\'es isol\'ees et de type A-D-E. 
\'Evidemment, $U(0)\subset U$. Nous construisons un champ alg\'ebrique de Deligne-Mumford
$\widetilde{U}:=\widetilde{U}(X,L)$  s\'epar\'e 
et propre sur son espace de module $U$ contenant $U(0)$ comme un ouvert de Zariski de sorte que,
le morphisme surjectif $\pi_1(U(0),\eta)\to \pi_1(\widetilde{U},\eta)$
(voir \cite{No1,No2}
pour les groupes d'homotopie des champs topologiques) associ\'e \`a 
l'inclusion de $U(0)$ dans $U$ a un noyau contenu dans celui de $\rho$. 
La repr\'esentation $\rho$ descend alors \`a une repr\'esentation 
$$\bar \rho:\pi_1(\widetilde{U}, \eta) \to O:=O (H_{ev}^n(X_{\eta}, \Q), \int_{X_{\eta}} - \cup -). $$
apparaissant comme l'holonomie d'une $\Q$-Variation de Structures de Hodge polaris\'ee de poids $n$
$(\V_{\widetilde{U}}, \mathcal{F}, S)$
sur $\widetilde{U}$. L'image $\Gamma$ de $\rho$ est la m\^eme que celle de $\bar \rho$ et par un th\'eor\`eme
classique de Beauville \cite{Bea}, c'est un sous-groupe arithm\'etique du groupe orthogonal $\mathbf{G}=O (H_{ev}^n(X_{\eta}, \R), \int_{X_{\eta}} - \cup -)$. 

On note par $\mathcal{D}:=\mathbf{G}/\mathbf{U}$ le domaine de Griffiths attach\'e \`a $(\V_{\widetilde{U}}, F, S)$ \cite{Gri},
$\mathbf{U}$ \'etant le sous-groupe qui stabilise la structure de Hodge 
sur $H_{ev}^n(X_{\eta}, \Q)$. On rappelle que $\mathcal{D}$ a une structure naturelle
 de vari\'et\'e complexe homog\`ene naturelle
et porte  une  distribution holomorphe horizontale $\mathbf{G}$-\'equivariante.
 L'action de $\Gamma$ sur $\mathcal{D}$ est proprement discontinue
et le champ quotient $[\Gamma \backslash \mathcal{D}]$ est un orbifold complexe. 
L'application des p\'eriodes de $(\V_{\widetilde{U}},  \mathcal{F}, S)$ d\'efinit une application holomorphe horizontale
de champs complexes analytiques $p:\widetilde{U} \to [\Gamma \backslash \mathcal{D}]$.

Pour toute vari\'et\'e $Z$ propre  sur $U$ (ce qui existe avec $\dim(Z)<c$), 
on peut construire un rev\^etement Galoisien fini $Z'\to Z$
de groupe $G$
\'etale au dessus de $U(1)\cap Z$ et un morphisme $[G\backslash Z']\to \widetilde{U}$ ce qui fournit
 un syst\`eme local $G$-\'equivariant $\mathbb{V}_{Z'}$ 
qui est sous jacent \`a  une $\Q$-VSH $G$-\'equivariante
sur $Z'$.  
Plus g\'en\'eralement pour tout morphisme $F:Y\to \widetilde{U} $ 
avec $Y$ projective-alg\'ebrique $\V_Y:=F^*(\V_{\widetilde{U}},  \mathcal{F}, S)$
est une $\Q$-VSH d'application de p\'eriodes $p\circ F$. 
Pour plus de d\'etails sur ces exemples,  diverses g\'en\'eralisations et une 
\'etude cohomologique de $\mathbb{V}_{\mathbb{P}'}$
si $\mathbb{P}\subset|L|$ est un espace lin\'eaire g\'en\'erique
de dimension $\le 6$ 
voir \cite{Megy}.

Soit $k\in \N$. 
L'ensemble $U(k)\subset |L|=\mathbb {P} H^0(X, L)$ des
 hypersurfaces $X_f:=\{f=0 \}$ 
\`a singularit\'es simples de nombre de Tjurina total
$\tau(f)$ inf\'erieur ou \'egal \`a $k$
est un ouvert dont le compl\'ementaire a codimension $c\ge\min(7,k)$
d\`es que $L$ est $k$-ample. On a $U(0)\subset U(k)\subset U$
et on note 
 $\wt{U(k)}=\widetilde{U}\times_U U(k)$. C'est un sous champ de $\widetilde{U}$ d'espace de modules  $U(k)$

  Le principal 
r\'esultat de cet article est un th\'eor\`eme de Torelli infinit\'esimal. Pour l'\'enoncer, 
nous avons besoin d'une d\'efinition. Si $\mathcal{I}\subset O_{X}$ est un id\'eal
coh\'erent, on note $p_X(\mathcal{I})=\min\{m, \ h^1(X,\mathcal{I}(m))=0\}$.
D'autre part, on note $H_k$ l'ensemble des id\'eaux coh\'erents dont le cosupport est un sous sch\'ema artinien 
de longueur $\le k$, et  
qui sont localement des id\'eaux jacobiens de singularit\'es
isol\'ees simples. Alors $\sup_{\mathcal{I}\in H_k} p_X(\mathcal{I})< \infty$. On note $s_k(X)$ ce supremum.

\begin{theoi}\label{torelli2} Si $X={\mathbb{P}}^{n+1}$, $n\ge 2$ pair, 
et $L=O_X(d)$ et $k\in \N$ v\'erifie $d\ge n+3+s_k(\mathbb{P}^{n+1})$,
 l'action de $PGL(n+2)$ sur $|L|$ se rel\`eve \`a $\widetilde{U}$ en pr\'eservant 
les $\wt{U(k)}$ et 
  la diff\'erentielle de la restriction de l'application de p\'eriodes
 de  $\mathbb{V}_{\wt{U}}$ \`a chaque strate  $\wt{U(k)}-\wt{U(k-1)}$ 
a pour noyau le tangent de l'orbite de $PGL(n+2)$.   
\end{theoi}

On a  $s_0(\mathbb{P}^{n+1})=s_1(\mathbb{P}^{n+1})=0$, $s_2(\mathbb{P}^{n+1})=2$. 
Nos  bornes ne sont pas optimales. Si les surfaces quintiques avec un n\oe ud sont
obtenues par notre th\'eor\`eme, ce dernier est vide pour les surfaces quartiques
\`a singularit\'es simples
alors que Torelli infinit\'esimal est bien connu dans ce cas. 
 L'obtention de bornes optimales n\'ecessiterait des 
arguments nettement plus fins non d\'evelopp\'es ici.

\begin{coroi}
Sous les hypoth\`eses pr\'ec\'edentes $p:[PGL(n+2)\backslash \wt{U(k)}]\to 
[\Gamma \backslash D]$ est finie.
\end{coroi}

Le r\'esultat avec $k=0$ est un r\'esultat classique de Griffiths \cite{Gri2}. 
La preuve du th\'eore\`eme \ref{torelli2} repose sur le calcul de la diff\'erentielle de l'application de p\'eriodes pour des 
hypersurfaces nodales issue du travail fondamental \cite{DS}
et de l'\'etude  de leur filtration de Hodge dans\cite{DSW}. 
Nous \'etendons une partie des r\'esultats de ces articles
aux hypersurfaces a singularit\'es simples. Cette extension effectu\'ee, l'\'enonc\'e de type
Torelli infinit\'esimal repose sur
 une variante donn\'ee au lemme \ref{lemac}
du th\'eor\`eme de Macaulay pour des hypersurfaces 
\`a singularit\'es isol\'ees quasi-homog\`enes,  exactement comme dans \cite{Voi}. 
Techniquement, nos r\'esultats sont compl\'ementaires de ceux de
 \cite{DSt1,DSt2} qui ne consid\`erent pas la question de Torelli infinit\'esimal. 
\footnote{Alors que nous finissions de r\'ediger ce travail, A. Dimca nous a
 signal\'e que
le lemme \ref{lemac} r\'esultait de \cite{DS2} qui traite le cas plus g\'en\'eral 
des singularit\'es isol\'ees quelconques, moyennant une traduction
 qui n'est pas si \'evidente pour nous. Notre preuve demandant moins de technologie
et restant assez courte, nous avons donc pr\'ef\'er\'e la conserver. }

Une g\'en\'eralisation du th\'eor\`eme de  Griffiths
sur les int\'egrales rationnelles donnant une interpr\'etation
de la diff\'erentielle de l'application de p\'eriodes comme op\'erateur de multiplication
pour les directions transverses aux strates isosinguli\`eres ne semble pas
avoir \'et\'e consid\'er\'ee de fa\c{c}on syst\'ematique dans la litt\'erature.  
De m\^eme il est probable que le th\'eor\`eme \ref{torelli2} se g\'en\'eralise 
pour $X$ quelconque pourvu que $L$ soit assez ample mais l\`a encore nous n'avons pas trouv\'e de r\'ef\'erence 
dans la litt\'erature. 
 Nous laissons ces questions pour de futures recherches.

L'application \`a la conjecture de Shafarevich est imm\'ediate en utilisant
une construction de \cite{Crelle, Eys}: 
\begin{coroi} \label{stein} Hypoth\`eses et notations comme au th\'eor\`eme \ref{torelli2}. 
 Soit $Z$ une vari\'et\'e projective lisse et $f: Z\to \wt{U(k+1)}$
un morphisme fini. Alors le rev\^etement universel de $Z$
est une vari\'et\'e de Stein.
\end{coroi}

Dans le cas o\`u $f$ est g\'en\'eriquement fini, l'\'etude de la 
conjecture de Shafarevich semble beaucoup plus d\'elicate et  nous ne savons pas non plus
la  d\'ecider dans tous les cas.

Nous tenons \`a remercier D. Barlet,  N. Borne, M. Brion, A. Dimca,  S. Druel,   L. Gruson, C. Peters, C. Voisin, M. H. Saito
et tout particuli\`erement A. Otwinowska pour 
d'utiles remarques sur les questions trait\'ees ici.

\section{Structure orbifold sur $U(k)$ et prolongement de la repr\'esentation de monodromie}

Soit $X$ une variété projective lisse complexe de dimension impaire $n+1$, $L$ un fibré en droites sans point base sur $X$, et $f\in H^0(X,L)-\{0\}$ 
telle que l'hypersurface $\{ f=0 \}:=X_f$ n'ait que des singularit\'es isol\'ees. 
Soit $\Sigma \subset X_f$ le sous sch\'ema artinien de $X$ de support ${X_f^{sing}}$ 
d\'efini par l'annulation du premier jet de $f$. La longueur de $\Sigma$, 
i.e. le nombre de Tjurina total $\tau(f)$,  est d\'efini par la relation 
$$\tau(f)=\sum_{P \in |\Sigma | } \tau_P (f)=
\sum_{P \in \Sigma} \dim(\OSp).$$
Choisissant des coordonn\'ees locales et une trivialisation locale de $L$ et notant  $f_P $ la fonction qui d\'efinit $f$
dans ces coordonn\'ees, 
 on voit que $\OSp$ est isomorphe \`a l'alg\`ebre de Tjurina
$O_{\C ^{n+1},0}/(f_P ,\frac{\partial f_P }{\partial x_1}, \ldots,\frac{\partial f_P }{\partial x_{n+1}})$. 

 Soit $k\in \N$, et $U(k)\subset |L|$  l'ouvert de Zariski constitu\'e des $[f]\in L$ tels que 
$X_f$ n'ait que des singularit\'es isol\'ees, simples (autrement dit de type A-D-E), et telles que  $\tau(f)\le k$.
On a $U(0)\subset U(1) \subset \ldots \subset U=\bigcup_kU(k) $. Remarquons que $U$ est non vide puisque l'ouvert $U(0)$ des $[f]$ tels que $X_f$ soit lisse est non vide par Bertini.\\

Dans la suite de cette section, on introduit (\ref{11}) un champ de Deligne-Mumford  muni d'une variation de structure de Hodge. Dans certains cas (\ref{13}), ce champ est également muni de l'action d'un groupe algébrique, et la variation de structure de Hodge descend au champ d'Artin quotient.

\subsection{Compl\'ements \`a  la construction de \cite{Megy}}
\label{11} Cette construction peut se reformuler et s'étendre comme suit :

\begin{prop} \label{rhobar}
Soient $X$,  $L$ et $U$ comme ci-dessus. Il existe un champ de Deligne-Mumford $\widetilde{U}= \bigcup_k \wt{U(k)}$
compactifiant $U(0)$ et  fini sur son espace de modules $U$,
tel que la repr\'esentation de monodromie de $\pi_1(U(0))$ se prolonge \`a 
$\bar\rho:=\bar{\rho}_{(X,L)}: \pi_1(\widetilde{U}, \eta)\to O(H^n(X_{\eta},\Q))$ o\`u $\eta\in U(0)$ est un point base arbitraire.
\end{prop}

\begin{prv}

\label{defmap}

Soit $f \in H^0(X,L)$ tel que l'hypersurface $X_f \subset X$ soit à singularités isolées. Tout voisinage de $[f]$ dans $\mathbb{P} H^0(X,L)$ induit une déformation de $X_f$ et donc, pour chaque singulier $P \in |\Sigma|$, une déformation du germe singulier $(X_f,P)$. On en déduit un morphisme de germes 
\[\lambda: (|L|,[f])\to \prod_{P \in |\Sigma|} Def(X_f,P )\]
o\`u $Def(X_f,P )$ est la base de \og la\fg{} d\'eformation miniverselle de la singularit\'e isol\'ee d'hypersurface $(X_f,P)$. Il est connu que dans ce cas, cette base $Def(X_f,P )$ est lisse, et naturellement isomorphe \`a un voisinage de $0$ dans
l'espace vectoriel $\OSp$. La d\'eformation $(X\times |L|, P \times [f])\supset (\mathcal{X}, P \times[f]) \to(|L|,[f])$ est induite
par la d\'eformation miniverselle \cite[Ch. 6]{Lo}. En particulier les \'el\'ements de $|L|$ proches de $[f]$ 
singuliers pr\`es de $P$ sont ceux que $\lambda$ envoie dans le discriminant de $Def(X_f,P )$. 
En choisissant une trivialisation de $L_P $ de fa\c{c}on ad\'equate, 
l'application $\lambda$
 a pour diff\'erentielle l'application d'évaluation des jets aux points singuliers de $X_f$ : 
\[
ev_{[f]}=\bigoplus_{P \in |\Sigma | } ev_P : T_{[f]} |L| \to \  \bigoplus_{P \in |\Sigma | } \OSp \otimes_{\OS} L
\]
qui provient par quotient de l'application naturelle $H^0(X,L) \to H^0(X,L\otimes \OS)$.\\

\emph{Stratification par le nombre de Tjurina.}\\

Soit $k\in \N$.  
Posons $Z(k)=U(k)-U(k-1)$ et supposons $[f]\in U(k)$. 
$Z(k)$ est localement ferm\'e dans  $U$. Pour $l\le k$ le germe $(Z(k),[f])$ se laisse
d\'ecrire comme l'image r\'eciproque par $\lambda$ de la r\'eunion sur toutes les partitions de $l$ 
$$l=\sum_{P \in |\Sigma|} l_P $$
de $\prod_{P \in |\Sigma|} Z_u(l_P )$ o\`u $Z_u(l_P )\subset Def(X_f,P )$ est le lieu des germes de nombre de Tjurina
 \'egal \`a $l_P $. 

Si $ev_{[f]}$ est surjectif, $\lambda$ est \'equivalent \`a la seconde projection du produit $(Z(k), [f])\times \prod_{P\in |\Sigma|} Def(X_f,P )$ 
\cite[Ch. 6]{Lo}
et la stratification $U(k)=\bigcup_{l=0}^k Z(l)$ est de Whithney en $[f]$ puisque
 c'est le cas pour la stratification par les nombres de Tjurina de la d\'eformation miniverselle d'un germe de singularit\'e simple.  
 La strate $Z(l)$ est alors de codimension $l$.\\

\emph{Rev\^etement Galoisien neutralisant la monodromie locale \`a l'infini\footnote{L'hypothèse $n\equiv 0[2]$ n'est utilisée qu'à partir de ce point.}}\\

Soit $B$ un voisinage  de $[f]$ suffisamment petit dans $|L|$.
Notons $\rho$
la repr\'esentation  de monodromie de $R^n p_{1*} \Q_{\mathcal{X}(1)}$ et $\rho_P $ 
la repr\'esentation  de $\pi_1(Def(X,P )-\overline{Z(1)})$ d\'efinie par la monodromie de la fibre de Milnor de $P\in|\Sigma|$. 
Alors,  $\rho|_{\pi_1(B-\overline{Z(1)})}$ se d\'ecompose comme somme directe d'un facteur trivial et
du tiré en arrière par $\lambda$ de $\oplus_P  \rho_P $. 
Puisque, en dimension paire, la monodromie locale $\rho_P $ est d'image finie, il suit que $\rho(\pi_1(B-\overline{Z(1)})$ est fini. 

Par un th\'eor\`eme de Selberg, il existe un sous-groupe normal sans torsion d'indice fini dans $\rho(\pi_1(|L|-\overline{Z(1)}))$. 
Le rev\^etement \'etale fini correspondant $\eta: U'(0)\to U-\overline{Z(1)}$
est galoisien de groupe $G$ et se prolonge par le th\'eor\`eme de Grauert-Remmert en un rev\^etement Galoisien
normal encore not\'e $\eta: |L|' \to |L|$. On note aussi $ U'(k)=|L|'\times_{|L|} U(k)$ et 
$\eta: U'(k)\to U(k)$ la restriction \`a cet ouvert. \\

\emph{Orbifold $\wt{U(k)}$.}\\

Consid\'erons le champ quotient $\wt{U(k)}= [U'(k)/G]$. 
Si $ev_{[f]}$ est surjective pour tout $[f]\in U(k)$,  
$U'(k)$ est lisse (voir \cite{Megy}) donc $\wt{U(k)}$ est un champ de Deligne-Mumford lisse
d'espace de modules $U(k)$ et l'isotropie en $[f]$
 est le produit des groupes de monodromie locale des singularit\'es. On note $\pi_1(\wt{U(k)})$ le  groupe fondamental du champ topologique sous-jacent \cite{No1,No2}.

Par construction, toute composante connexe $V$ de la pr\'eimage  de $U-Z(1)$ dans $U'(k)$, v\'erifie $\rho(\pi_1(V))=\{e \}$
et que $U'(k)$ soit  lisse ou non, la repr\'esentation $\rho \circ \eta_*$ est induite par une repr\'esentation
$\rho':\pi_1(U'(k),x)\to O(H^n(X_g,\Q))$ o\`u $x\in U'(1)$ est un point base et $[g]=\eta(x)$. 
Le syst\`eme local correspondant est $G$-\'equivariant et descend \`a un syst\`eme
local  $\V_{\wt{U(k)}}$ sur $\wt{U(k)}$ et on notera $\bar\rho: \pi_1(\wt{U(k)}, x)\to O(H^n(X_g,\Q))$
sa monodromie.
Si $\wt{U(k)}$ est lisse ce syst\`eme local est
sous jacent \`a une $\Q$-VSH  par \cite{Gri}.
\end{prv}
\medskip

Notre discussion implique aussi la pr\'ecision suivante sauf pour le dernier point pour lequel
on renvoie \`a \cite{Megy}. 

\begin{prop}
Si, de plus,  $L$ est $k$-jet-ample, $\wt{U(k)}$ est lisse, $\mathrm{codim}_{|L|} U(k)^c\ge \max(7,k)$ et la stratification $\wt{U(k)}=\cup_{l=0}^k \wt{Z(l)}$
(o\`u  l'on a pos\'e $\wt{Z(l)}=\wt{U(l)}-\wt{U(l+1)}$)
est de Whithney. 
\end{prop}

Notons que $\gamma:\wt{Z(k)} \to Z(k)$ est une gerbe de lien un 
sch\'ema en groupes fini si,  pour tout  $[f]\in Z(k)$,   $ev_{[f]}$ est surjectif. 
Ceci permet de d\'efinir un syst\`eme local $\V_{Z(k)}:=\gamma_* \V_{\wt{Z(k)}}$ qui est 
sous jacent \`a une VSH polaris\'ee sur $Z(k)$.

\subsection{Une question ouverte}\label{12}

Dans le cas $X=\mathbb{P}^3$,  $L=O_{\mathbb{P}^3}(4)$ on sait que $\bar{\rho}$ est un isomorphisme sur son image
gr\^ace au th\'eor\`eme de Torelli pour les surfaces K3. Une conjecture de Carlson-Toledo
pr\'edit que les seuls r\'eseaux de groupes alg\'ebriques r\'eels semisimples
apparaissant comme groupes k\"ahl\'eriens sont ceux des groupes de type hermitiens sym\'etriques. 
Cette conjecture implique que $\bar{\rho}$ n'est pas injective pour $d\ge 5$. 
Ceci motive:
\begin{conj}
Pour $d\ge 5$ le noyau de $\bar{\rho}=\bar{\rho}_{(\mathbb{P}^3,O_{\mathbb{P}^3}(d))}$ est un groupe
infini.  
\end{conj}
La m\'ethode de \cite{CTDuke} pour prouver le fait analogue dans le cas de $\rho$
ne s'applique malheureusement pas ici. Nous ne voyons pas comment construire d'autres
repr\'esentations lin\'eaires de $\pi_1(\widetilde{U}(\mathbb{P}^3,O_{\mathbb{P}^3}(d)),\eta)$. 
Ceci motive la:
\begin{ques}
 Le  groupe  $\pi_1(\widetilde{U}(\mathbb{P}^3,O_{\mathbb{P}^3}(d)),\eta)$ admet il d'autres repr\'esentations
complexes 
que les repr\'esentations de la forme $\alpha\circ \bar \rho$ o\`u $\alpha$ est une repr\'esentation rationnelle
de $O(H^n(X_g,\R))$?
\end{ques}
Pour $d=4$ la r\'eponse \`a cette question est n\'egative par le th\'eor\`eme de superrigidit\'e de Margulis.

\subsection{Rel\`evement de l'action de $PAut(X,L)$}\label{13}

Le groupe $Aut(X,L)$ des automorphismes du couple $(X,L)$ agit sur $|L|$ \`a travers $PAut(X,L)=Aut(X,L)/\C^*$ 
 en pr\'eservant les $U(k)$.

\begin{prop}
 Si le groupe $PAut(X,L)$ est semisimple, son action sur $U$ se rel\`eve \`a une action sur $\wt{U}$. 
\end{prop}
\begin{prv}
 L'alg\`ebre de Lie ${\mathfrak{paut}}(X,L)^{op} \subset  H^0(\Theta_{|L|})$ 
est une alg\`ebre de Lie de champs de vecteurs holomorphes qui sont  tangents \`a chaque  $Z(k)$.
Or $\overline{Z(1)}$ est le lieu de ramification de $|L|'\to |L|$. 
Donc ces champs de vecteurs se rel\`event \`a des champs de vecteurs sur $|L|'$ car tout champ de vecteurs
tangent au lieu de ramification de $Z\to Y$ avec $Y$ lisse et $Z$ normal se rel\`eve \`a $Z$.

En effet,  le rel\`evement a lieu en codimension un car un germe de champ de vecteurs de la forme $a(z,(w_j))z\frac{\partial}{\partial z}+ b_i(z,(w_j))\frac{\partial}{\partial w_i}$
$a,b\in \C\{z,w_1, \ldots,w_n\}$ se rel\`eve bien \`a un germe de champ de vecteurs holomorphes par
 un morphisme de la forme $(z,(w_j)\to (z^e,(w_j))$. Ce rel\`evement en codimension un se prolonge à $Y$ tout entier car le faisceau des champs de vecteurs holomorphes sur l'espace normal $Y$ est r\'eflexif.

On dispose donc d'un morphisme d'alg\`ebres de Lie ${\mathfrak{paut}}(X,L)^{op} \to  H^0(\Theta_{U'})^G$ qui s'exponentie
en un morphisme de groupes de Lie complexes 
du rev\^etement universel topologique $\tilde{P}$ de $PAut(X,L)$
vers le centralisateur $C(G, Aut(U')$ de  $G$ dans $Aut(U')$.  

Si $PAut(X,L)$ est semisimple, le groupe $\wt P $ est un groupe alg\'ebrique affine semisimple 
et le morphisme correspondant est un morphisme de groupes alg\'ebriques. 
Le noyau $N$ du morphisme $\wt P \to Aut'(U)$ est contenu dans le noyau $N'$ de $\wt P \to  PAut(X,L)$
car l'action de $\wt P$ redescend \`a une action de $\wt P$ sur $U$ factorisant via $PAut(X,L)$. 
Mais $N'/N$ commute \`a $G$ et pr\'eserve $\eta$. Comme $\eta$ est Galoisien de groupe $G$
on d\'eduit que $N'/N\subset Z(G)$ et que donc, en divisant $U'$ par $N'/N$, 
on obtient une action de $PAut(X,L)$ sur $U'':=U'\slash (N'/N)$ qui commute \`a $Gal(U''/U)=G\slash(N'/N):=G'$
et descend \`a une action de $PAut(X,L)$ sur $\wt{U}\simeq [U'' \slash G']$. 
\end{prv}

\begin{coro}
 L'application des p\'eriodes  $\mathcal P:\widetilde{U} \to [\Gamma \backslash \mathcal{D}]$ attach\'ee \`a $\bar \rho$ descend \`a une application d\'efinie
sur le champ quotient
$[PAut(X,L) \backslash \widetilde{U}]\to [\Gamma \backslash \mathcal{D}]$.
\end{coro}

\section{Interpr\'etation g\'eom\'etrique de la repr\'esentation de monodromie prolong\'ee}
\label{sec2}

Avec les notations de \ref{11}, soit $[f]\in U$ et $X_f\subset X$ l'hypersurface de
 $X$ \`a singularit\'es isol\'ees simples d\'efinie par $f$ et $ev_{[f]}$ la flèche d'évaluation des jets aux points singuliers de $X_f$.
 
 \begin{prop}\label{prop21} Supposons $ev_{[f]}$ surjectif pour tout $X_f\in Z(k)$. 
Les sous faisceaux de $\Theta_{Z(k)}$ d\'efinis par les noyaux des $Gr_F (\nabla)$ pour les  VSH sous-jacentes aux syst\`emes locaux
$Gr^{n+2}_W R^{n+1} (p_1)_* \Q_{(Z(k)\times  X- \mathcal{X}_{Z(k)})}$ et $\V_{Z(k)}$
sont \'egaux.
\end{prop}

Remarquons que l'hypothèse de la proposition est satisfaite si par exemple $L$ est le produit tensoriel d'au moins $\tau(f)$ fibrés tr\`es amples.

Pour montrer la proposition, on compare la cohomologie de $X_f$ à celle de son compl\'ementaire  (\ref{21}) ainsi qu'à la structure limite (\ref{22}), d'abord ponctuellement puis en famille (\ref{23}), ce qui donne le résultat. Une preuve alternative et plus constructive est donnée dans le cas de familles de surfaces en (\ref{24}).
 
\subsection{Relation entre les structures de Hodge de $X_f$ et de $X-X_f$}\label{21}

Une singularit\'e simple a une forme d'intersection d\'efinie n\'egative, en particulier non d\'eg\'en\'er\'ee, 
donc $X_f$ est une vari\'et\'e d'homologie rationnelle \cite[prop. 4.7]{Dimca92}. 
Le groupe $H^n(X_f,\Q)$ co\"incide donc avec le groupe de cohomologie d'intersection $IH^n(X_f,\Q)$ et porte une structure de Hodge pure. 
Le groupe $H^{n+1}(X-X_f, \Z)$ 
 porte une structure de Hodge mixte  \cite{Del} de poids $n+1$ et $n+2$, 
entrant dans une suite exacte \label{sec0}:
\[H^{n+1}(X, \Z) \to H^{n+1}(X-X_f, \Z)\to H^{n+2}_{X_f}(X,\Z) \to H^{n+2}(X,\Z) \leqno{(\ref{sec0})}
\]
Toujours parce que $X_f$ est une vari\'et\'e d'homologie rationnelle, on a un isomorphisme de structure de Hodge $H^n(X_f, \Q) (1)\to H^{n+2}_{X_f}(X,\Q)$. 
Ainsi $Gr_W^{n+2} H^{n+1}(X-X_f, \Z)$ est une sous structure de Hodge pure de poids $n+2$
de $H^{n+2}_{X_f}(X,\Z)$ et sa filtration de Hodge v\'erifie 
\[0=F^{n+2} \subset F^{n+1} \subset \ldots \subset F^1= F^0=Gr_W^{n+2} H^{n+1}(X-X_f, \Z)
\]

On note que si $H^{n+1}(X, \Z)=0$, ce qui est le cas si $X=\mathbb{P}^{n+1}$ puisque $n$ est pair, alors $H^{n+1}(X-X_f, \Z)$
est pure de poids $n+2$. 

\subsection{Structure de Hodge de $X_f$ et  structure de Hodge limite}\label{22}

Soit $i: \Delta \to U$ un disque analytique  tel que $i(0)=[f]$ et $i(\Delta^*)\subset U(0)$. On pose $\mathcal{X}_{\Delta}=\mathcal{X}\times_U \Delta$.
Alors $t:\mathcal{X}_{\Delta}\to \Delta$ est un morphisme projectif plat lisse hors de $\{t=0\}$
et $(\mathcal{X}_{\Delta})_0=X_f$. Si la monodromie de $\mathcal{X}_{\Delta}\to \Delta^*$ est unipotente, c'est \`a dire triviale puisque les groupes
de monodromie locale pr\`es de $[f]$ sont finis, ceci  permet de d\'efinir la structure de Hodge limite $H^n_{lim}(X_f, \Q)$
au sens de \cite{St}.

Gr\^ace \`a \cite{MHM}, 
on a un isomorphisme de structures de Hodge pures
\[H^n_{lim}(X_f, \Q)=\mathbb{H}^0(X_f, \psi_t(\Q_{\mathcal{X}_{\Delta}}[n]))\] 

L'objet $\phi_t(\Q_{\mathcal{X}_{\Delta}}[n])$ est concentr\'e aux points singuliers de $X_f$. 
Le Module de Hodge Mixte $\psi_t(\Q_{\mathcal{X}_{\Delta}}[n])=\psi_t(IC_{\mathcal{X}_{\Delta}}(\Q))$
est un module de Hodge polarisable pur car $IC_{\mathcal{X}_{\Delta}}(\Q)$ est un module de
Hodge polarisable et que le logarithme de la monodromie v\'erifie $N=0$, en vertu de \cite[(0.7), p. 852]{MHP}.

La premi\`ere fl\`eche $a$ du   triangle distingu\'e canonique 
\[  \Q_{X_f}[n] \to \psi_t(\Q_{\mathcal{X}_{\Delta}}) \xrightarrow{can} \phi_t(\Q_{\mathcal{X}_{\Delta}})\xrightarrow{+1}\]
est donc un morphisme de Modules de Hodge Polarisables de m\^eme poids. 
La cat\'egorie des modules de Hodge 
Polarisables de poids donn\'e \'etant ab\'elienne et semi-simple \cite[Lemme 5, p. 854]{MHP}
il suit que $\psi_t(\Q_{\mathcal{X}_{\Delta}})\simeq \Q_{X_f}[n] \oplus \mathrm{Coker}(a)$ o\`u
$\mathrm{Coker}(a)$ a m\^eme poids que $\Q_{X_f}[n]$. Il suit que $\phi_t(\Q_{\mathcal{X}_{\Delta}})$
est isomorphe \`a $\mathrm{Coker}(a)$. 

Ceci fournit une suite exacte scind\'ee de structures de Hodge pures \label{sec}:
\[ 0\to H^n (X_f, \Q)  \to H^n_{lim}(X_f, \Q) \to \bigoplus_{x\in |\Sigma|} 
\mathcal H^n(\phi_t(\Q_{\mathcal{X}_{\Delta}}))_x\to 0 \leqno{(\ref{sec})}\]

\subsection{Comparaison avec la variation de structure de Hodge des vari\'et\'es singuli\`eres} \label{23}

Soit avec les notations du paragraphe \ref{defmap} 
un germe de disque analytique $j: (\Delta,0) \to \prod_{p\in|\Sigma|} Def(X_f,p)$
tel que $j(\Delta^*)\subset \prod_{p\in|\Sigma|} Def(X_f,p)(0)$ et
dont la monodromie est nulle. 

On d\'efinit un germe $ T:=(|L|,[f])\times_{\lambda, j} \Delta$
muni de la famille d'hypersurfaces $X\times T\supset \mathcal{X}_T\buildrel{\pi}\over{\to} T$ .

La projection naturelle $t:T\to \Delta$ d\'efinit une fonction holomorphe et $S:=\{t=0\}$ est exactement le germe
en $[f]$ de la strate isosinguli\`ere $Z(k)$ de $X_f$. Notons enfin $g$ la composée $\mathcal{X}_T \xrightarrow{\pi} T \xrightarrow{t}\Delta$.

On a sur $\mathcal{X}_S$ un triangle distingué

\[
  \Q_{\mathcal{X}_S}[n] \xrightarrow{sp} \psi_g(\Q_{\mathcal{X}_T}) \xrightarrow{can}\phi_g(\Q_{\mathcal{X}_T})\xrightarrow{+1} .\]

Le dernier terme est supporté sur le lieu  singulier de $g$ qui est une union finie de 
multisections \'etales de $\mathcal{X}_S \to S$
qui d\'ecrivent les points singuliers des $X_s$ pour $s\in S$. 
Comme ces singularités sont de type ADE donc rigides, 
le dernier terme est donc une somme directe, 
  de tir\'es en arrière de faisceaux gratte-ciel sur la d\'eformation universelle locale de chaque singularit\'e.

En appliquant $\pi_*$, et en utilisant l' isomorphisme naturel $\pi_*\psi_g \simeq \psi_t \pi_*$,  on obtient sur $S$  le triangle
\[
 \pi_*\Q_{\mathcal{X}_S} \to \psi_t(\pi_*\Q_{\mathcal{X}_T})\to \pi_*\phi_g(\Q_{\mathcal{X}_T}) \xrightarrow{+1}  .\]
Ce dernier donne en cohomologie une suite exacte  sur $S$ 
\[ 0\to R^n\pi_*\Q_{\mathcal{X}_S} \to  \V_S\oplus H^n(X,\Q) \to  \bigoplus_{x\in |\Sigma|} \mathcal R^n\pi_*\phi_g(\Q_{\mathcal{X}_{\Delta}})_x.\]

La suite exacte (\ref{sec}) implique que les tiges des termes de cette suite exacte sont
des structures de Hodge pures de m\^eme poids et donc que cette suite exacte est une suite
exacte de Modules de Hodge polarisables purs qui est scind\'ee toujours par  \cite[Lemme5, p. 854]{MHP}. 

Finalement, en tout point de $S$ le noyau de $Gr_F (\nabla)$
 pour la $\Q$-Variation de Structures de Hodge $R^n\pi_*\Q_{\mathcal{X}_S}$ et celui de  $\V_S$
sont les m\^emes. 

Compte tenu de la suite exacte (\ref{sec0}), ceci \'etablit la proposition \ref{prop21}.\qed

\subsection{Interpr\'etation dans le cas $n=2$} \label{24}

Donnons un argument alternatif permettant de d\'emontrer la proposition \ref{prop21} dans le cas $n=2$. 

Si on applique le th\'eor\`eme de r\'esolution simultan\'ee des singularit\'es DuVal 
\cite{Art} \cite{Bris3}, cf. \cite[p. 135]{KM} \`a la famille universelle $p_1:\mathcal{X}\to U$
on trouve un rev\^etement ramifi\'e $r:U^*\to U$ et
une application holomorphe propre et lisse $\bar u:\mathcal{X}^*\to U^*$ avec un morphisme 
$\pi:\mathcal{X}^*\to \mathcal{X}
\times_{U}U^*$ qui soit une r\'esolution simultan\'ee, c'est \`a dire que
pour tout $s\in U^*$ $\pi_s: \mathcal{X}^*_s\to \mathcal{X}_s$ soit une r\'esolution. 
De plus, ces r\'esolutions peuvent \^etre suppos\'ees minimales. 
Le morphisme $r$ factorise par un morphisme fini $r':U^*\to \widetilde{U}$. 
 L'existence globale de  $r$ n'est pas \'evidente 
et il n'est pas clair pour nous que $r$ puisse \^etre choisi de fa\c{c}on \`a ce que $r'$ soit \'etale.
Toutefois, cf \cite{BW}, 
c'est le cas si on restreint $u$ \`a un petit voisinage d'un point de $U$.
Ceci implique que $\V_{\widetilde{U},[f]} \simeq H^2(X'_{f}, \Q)$ comme structures de Hodge
o\`u $X'_f\to X_f$ est la r\'esolution minimale
mais aussi que $r'^*\V_{\widetilde{U}}\simeq R^2\bar u_{*}\Q_{\mathcal{X}^*}$ comme variations de structures
de Hodge polarisables. 

Notons $Z^*(k)=r^{-1}( Z(k)-Z(k+1))$ et appelons $\mathcal{X}^*(k)$, resp. $\mathcal{X}(k)$ la restriction de 
$\mathcal{X}^*$, resp. $\mathcal{X}$ \`a $Z^*(k)$. $\mathcal{X}^*(k)$ est lisse sur $Z(k)$
 et, quitte \`a faire un rev\^etement \'etale de $Z^0(k)$, on peut supposer
 que  sur chaque composante connexe $Z^0(k)$ de $Z^*(k)$ 
l'ensemble singulier de $\mathcal{X}(k)$ est un produit de $Z^0(k)$ par un ensemble fini.
L'ensemble exceptionnel de
 $\pi_k:\mathcal{X}^*(k)\to\mathcal{X}(k)$ est globalement
 un produit de $Z^0(k)$ par une r\'eunion de configurations de courbes rationnelles du type A-D-E ad\'equat. 

Ceci donne des suites exactes de Variation de Structure de Hodge:
$$
0\to \Q(1)^{\oplus k}_{Z^0(k)} \to R^2\bar u_* \Q_{\mathcal{X}^*(k)}|_{Z^0(k)} \simeq(r')^*\V_{\widetilde{U}} |_{Z^0(k)}
\to R^2(p_1)_* \Q_{\mathcal{X}(k)}  \to 0
$$
Notons $v: U^*\times X - \mathcal{X}^* \to U^*$ la premi\`ere projection. La suite exacte (\ref{sec0}) implique que 
$R^2(p_1)_* \Q_{\mathcal{X}(k)}(1)\simeq Gr^4_W R^3 v_* \Q_{(U^*\times X - \mathcal{X}^*)}|_{Z^0(k)}$, 
puis par le th\'eor\`eme de semisimplicit\'e:
$$ (r')^*\V_{\widetilde{U}} |_{Z^0(k)} \simeq  \Q(1)^{\oplus k}_{Z^0(k)} \oplus
Gr^4_W R^3 v_* \Q_{(U^*\times X - \mathcal{X}^*)}|_{Z^0(k)}
$$
Ceci implique que sur $Z^0(k)$, les VSH $Gr^4_W R^3 v_* \Q_{(U^*\times X - \mathcal{X})}|_{Z^0(k)}$ et
 $\V_{\widetilde{U}}|_{Z^0(k)}$
ne diff\'erent que par un syst\`eme local de monodromie finie donc d'application de p\'eriodes constante.
 En particulier, Torelli infinit\'esimal pour $ \V_{\widetilde{U}} |_{Z^0(k)}$ \'equivaut \`a Torelli infinit\'esimal 
pour $Gr^4_W R^3 v_* \Q_{X\times U^* - \mathcal{X}}|_{Z^0(k)}$ ce qui \'equivaut \`a la proposition \ref{prop21}
dans  ce cas.

\section{Th\'eor\`eme de Macaulay avec singularit\'es mod\'er\'ees}
\label{3}

Claire Voisin nous a signal\'e qu'une variante du th\'eor\`eme de Macaulay 
autorisant un peu de singularit\'es devrait suivre en adaptant  \cite[pp. 427-428]{Voi}.  
Mettons en oeuvre de cette suggestion : apr\`es quelques g\'en\'eralit\'es sur le complexe de Koszul,
 nous obtenons la propri\'et\'e d'injectivit\'e \ref{injectivite1} en corollaire du r\'esultat 
de dualit\'e \ref{lemac}. Cela est suffisant pour montrer l'injectivit\'e de
 l'application des p\'eriodes de $\V_{\wt{U}}$ sur certaines strates de $\wt{U}$.

\subsection{Un lemme sur le complexe de Koszul d'une presqu'intersection compl\`ete}\label{31}

Soit $D\in\N^*$ un entier positif. Soit $G\subset H^0(X, L^D)$ un sous espace de dimension $n+2$. 
Notons $\mathcal{G}\subset \OX$ le faisceau d'id\'eaux engendr\'e par $G$, et $\Sigma\subset X$ le sous sch\'ema  tel que $\OS=\OX/\mathcal{G}$.

Posons $A=A(X,L)=\oplus_{k\in \N} H^0(X, L^k)$. 
On appelle $I$ l'id\'eal gradu\'e de $A(X,L)$ d\'efini en degr\'e $k$ par $I_k=H^0(\mathcal{G}(k))$. On a $(G)\subset I$.

Soit $K(G)_m^{\hdot}$ le complexe de Koszul en degr\'e $m$:
\begin{equation*}
H^0(L^{m -(n+2)D}) \otimes \Lambda^{n+2} G\to \ldots \to H^0(L^{m-D})\otimes G \to H^0(L^m)
\end{equation*}
o\`u le premier terme du complexe est par convention en degr\'e $0$. 
Consid\'erons \'egalement, avec la m\^eme convention,  sa version faisceautique $\mathcal{K}(G)_m^{\hdot}$:
\begin{equation*}
\OX(L^{m -(n+2)D}) \otimes \Lambda^{n+2} G\to \ldots \to \OX(L^{m-D})\otimes G \to \OX(L^m)
\end{equation*}

L'hypercohomologie de ce complexe de faisceaux est d\'ecrite par le lemme suivant.

\begin{lem} \label{le31}Si $\Sigma\subset X$ est artinien et 
localement d'intersection compl\`ete, alors, pour tout $m\in \Z$,  
$\mathbb{H}^i(\mathcal{K}(G)_m^{\hdot})=0$ pour $i\not= n+2,n+1$ et 
\begin{equation}
 \mathbb{H}^{n+2}(\mathcal{K}(G)_m^{\hdot})\simeq\mathbb{H}^{n+1}(\mathcal{K}(G)_m^{\hdot}) \simeq \bigoplus_{P \in |\Sigma|} \OSp.
\end{equation}
\end{lem}

\begin{prv} Le support $|\Sigma|$ de $\Sigma$
consiste en un nombre fini de points de $X$ et 
$$\OS= \bigoplus_{P \in |\Sigma|} \OSp$$ est une somme de faisceaux gratte-ciel. 
Par abus de langage on identifie $\OSp$ et l'alg\`ebre artinienne locale de $\Sigma$ en $P $. 
Ensuite, pour tout $P \in |\Sigma|$, 
d\'esignant par $\mathfrak{m}_P \subset \OXp$ l'id\'eal maximal, on a
$\dim_{\C} \G_P  \slash \mathfrak{m}_P \G_P = n+1$ et 
toute famille $(g_0, \ldots g_n)$ dans $\G_P $ induisant une base de $\G_P  \slash \mathfrak{m}_P \G_P $
 est une suite r\'eguli\`ere dans $\OXp$ engendrant $\G_P $
voir par exemple \cite[Thm 129]{Kap}. On peut donc choisir une base  $(g_1, \ldots, g_n, g_{n+1})$ de $G$ de sorte que 
  $(g_1, \ldots g_n)$ est une suite r\'eguli\`ere engendrant $\G_P $. 
Dans la cat\'egorie d\'eriv\'ee $D^b(Mod(\OXp))$ on a donc 
\begin{eqnarray*}                                                          
\mathcal{K}(G)_{0,P }^{\hdot} &=&
 \mathcal{K}(g_{1,P }, \ldots, g_{n+1,P })_0^{\hdot} \otimes  \mathcal{K}(g_{n+2,P })_0^{\hdot}\\
&\simeq& \OSp [-n-1] \otimes^L \mathcal{K}(g_{n+2,P })^{\hdot} \\ &\simeq& (\OSp[-n-1]\buildrel{0}\over\to \OSp [-n-2])
 \end{eqnarray*}
puis
$H^ {n+1}(\mathcal{K}(G)_{0,P }^{\hdot})\simeq H^{n+2}(\mathcal{K}(G)_{0,P }^{\hdot})\simeq \OSp$ comme $\OXp$-modules, les autres faisceaux de cohomologie \'etant
nuls. Par suite, les faisceaux de cohomologie non nuls de $\mathcal{K}(G)_m^{\hdot}$ sont les gratte-ciels
$H^{n+1}(\mathcal{K}(G)_m^{\hdot})\simeq H^{n+2}(\mathcal{K}(G)_m^{\hdot}) \simeq  \bigoplus_{P\in |\Sigma|} \OSp$
puisque $\mathcal{K}(G)_0^{\hdot}$ est acyclique hors de $|\Sigma|$. La suite spectrale d'hypercohomologie de $\mathcal{K}(G)_m^{\hdot}$
n'a donc qu'un seul terme non nul en $E_1$ qui est $d_1: H^0(H^{n+1}(\mathcal{K}(G)_m^{\hdot}))\to H^0(H^{n+2}(\mathcal{K}(G)_m^{\hdot}))$
et $d_1=0$ car c'est le cas apr\`es localisation.
Donc, pour tout $m\in \Z$,  $\mathbb{H}^i(\mathcal{K}(G)_m^{\hdot})=0$ pour $i\not= n+2,n+1$ et 
\begin{equation}\label{hyper}
 \mathbb{H}^{n+2}(\mathcal{K}(G)_m^{\hdot})\simeq\mathbb{H}^{n+1}(\mathcal{K}(G)_m^{\hdot}) \simeq \bigoplus_{P \in |\Sigma|} \OSp.
\end{equation}

\end{prv}

\label{suite_spectrale_bete}Par ailleurs, on a la filtration b\^ete de $\mathcal{K}(G)_0$ d\'efinie par 
$$\sigma_{\ge p}  \mathcal{K}(G)^{\hdot} = \OX(L^{-pD}) \otimes \Lambda^{n+1-p} G \to \ldots \to \OX.
$$
Cette filtration d\'ecroissante induit une filtration sur $\mathbb{H}^*(\mathcal{K}(G)_m^{\hdot})$ et la suite spectrale
correspondante a un terme $E_1^{p,q}= H^q(Gr_\sigma^p(\mathcal{K}(G)^{\hdot}))$. 
Elle d\'eg\'en\`ere en $E_{n+3}$. On remarque \'egalement que  $(E_1^{\hdot,0},d_1)=K(G)^{\hdot}_m$, 
de sorte que $E_2^{p,0}= H^p(K(G)^{\hdot}_m)$. Notamment, $E_2^{n+2,0}= H^{n+2}(K(G)^{\hdot}_m)= 
(R/(G))_m$ et une variation l\'eg\`ere de la preuve du lemme \ref{le31} permet de voir que
$G$ a une suite r\'eguli\`ere de longueur $n+1$ et donc que $E_2^{p,0}=H^{p}(K(G)^{\hdot}_m)=0$ pour $p\leq n$.

\subsection{Dualit\'e de Macaulay pour des hypersurfaces de $\PP$ \`a singularit\'es isol\'ees quasihomog\`enes}\label{32}

On se place dans le cas particulier suivant. Soit $f\in H^0(\PP,\OP(d))-\{0\}$ telle que l'hypersurface $X_f:=\{ f=0 \}$ n'aie que des singularit\'es isol\'ees quasi-homog\`enes. 
On consid\`ere $G\subset H^0(\PP,\OP(d-1))$ le sous espace vectoriel engendr\'e par les d\'eriv\'ees partielles de  $f$, et $(G)=J$ est alors l'id\'eal jacobien de $f$. 
Avec les notations pr\'ec\'edentes, $D = d-1$, et $\Sigma \subset X_f$
 s'identifie au sous sch\'ema artinien de $X$ de support ${X_f^{sing}}$ 
d\'efini par l'annulation du premier jet de $f$, $I$ \`a l'id\'eal de $A$
s'annulant sur $\Sigma$ et
 $\G=\mathcal{J}_{\Sigma}$ s'identifie, apr\`es introduction de coordonn\'ees locales, 
\`a l'id\'eal de Tjurina de la fonction correspondant \`a $f$. Puisque les singularit\'es de $X_f$ 
sont  quasi-homog\`enes, les id\'eaux de Milnor et de Tjurina co\"{i}ncident donc $\Sigma$ est 
intersection compl\`ete locale.
 
Si $l \in \Z$, on note 
\[\text{ev}_l : H^0(\PP,\OP(l)) \to H^0(\PP,\OS(l))\]
la fl\'eche induite par la surjection de faisceaux $\OP\to \OS$.

On note enfin $\sigma = (n+2)(d-2)$.
 
Sur l'espace projectif $\PP$, les faisceaux inversibles n'ont de cohomologie qu'en degr\'e $0$ ou $n+1$, on en d\'eduit l'annulation de la plupart des termes de la suite spectrale introduite en \ref{suite_spectrale_bete} : on a  $E_r^{p,q}=0$ sauf si $q=0$ ou $q=n+1$, et donc les seules fl\'eches $d_r$ non nulles sont les fl\'eches $d_1$ et la fl\'eche $d_{n+2} : E_{n+2}^{0,n+1}\to E_{n+2}^{n+2,0}$. Les deux propositions suivantes pr\'ecisent son image et sa coimage.

\begin{lem} \label{lemac}
La fl\'eche $d_{n+2} : E_{n+2}^{0,n+1}\to E_{n+2}^{n+2,0}$ induit un isomorphisme entre $\left(I\slash (G)\right)_{\sigma-m}^\vee$ et $\left(I\slash (G)\right)_m$.
\end{lem}

\begin{prv}
Cette fl\'eche induit un isomorphisme entre sa coimage et son image. Il suffit de calculer ces deux espaces.\\
\\
\emph{Premi\`ere \'etape: l'image de $d_{n+2}$ dans $E_{n+2}^{n+2,0}=R\slash(G)_m$ est  $\left(I\slash (G)\right)_m$.}\\
\\
On a 
 \begin{align*}
 {\rm Im}d_{n+2} &= \ker \left(E^{n+2,0}_{n+2}\twoheadrightarrow  E^{n+2,0}_{n+3}\right)\\
 &= \frac{\ker \left(E^{n+2,0}_{1}\twoheadrightarrow  E^{n+2,0}_{n+3}\right)}{\ker \left(E^{n+2,0}_{1}\twoheadrightarrow  E^{n+2,0}_{n+2}\right)}\\
 &= \frac{\ker \left(E^{n+2,0}_{1}\twoheadrightarrow  E^{n+2,0}_{n+3} \hookrightarrow \mathbb{H}^{n+2}(\mathcal{K}(G)_m^{\hdot})\right)}{\ker \left(E^{n+2,0}_{1}\twoheadrightarrow  E^{n+2,0}_{n+2}\right)}
 \end{align*}

Or, la fl\`eche compos\'ee 
\[ E^{n+2,0}_1 \twoheadrightarrow E^{n+2,0}_2 =  E^{n+2,0}_{n+2}\twoheadrightarrow  E^{n+2,0}_{n+3} \hookrightarrow \mathbb{H}^{n+2}(\mathcal{K}(G)_m^{\hdot})\]
apparaissant au num\'erateur co\"incide avec $H^0(\mathcal O(m)) \to H^0(\mathcal O / \mathcal G(m))$. Son noyau est donc $H^0(\mathcal G(m))=I_m$. D'autre part, au d\'enominateur, on a 
\[\ker \left(E^{n+2,0}_{1}\twoheadrightarrow  E^{n+2,0}_{n+2}\right) = \ker \left(E^{n+2,0}_{1}\twoheadrightarrow  E^{n+2,0}_{2}\right) = {\rm Im } d_1.\]
On en d\'eduit
\[
 {\rm Im}\:d_{n+2}
  = \frac{I_m}{{\rm Im }\: d_1}=\frac{I_m}{(G)_m}.\]

\emph{Deuxi\`eme \'etape : la coimage de $d_{n+2} : E_{n+2}^{0,n+1} \to E_{n+2}^{n+2,0}$ 
s'identifie par dualit\'e de Serre 
\`a $\left(I\slash (G)\right)_{-m+(n+2)(d-2)}^\vee = \left(I\slash (G)\right)_{\sigma-m}^\vee$.}

 La suite spectrale
duale $((E^{p,q}_r)^{\vee}, ^td_r)$
 s'identifie par dualit\'e de Serre \`a une renum\'erotation de
la suite spectrale de la filtration b\^ete de:
\begin{equation*}
\OmX(L^{-m}) \to  \OmX(L^{-m+d-1})\otimes G^{\vee} \to  \ldots \to \OmX(L^{-m +(n+2)(d-1)}) \otimes \Lambda^{n+2} G^{\vee}
\end{equation*}
Le complexe obtenu en tensorisant  par la droite complexe $\Lambda^{n+2} G\simeq \C $
s'identifie \`a $\mathcal{K}(G)_{-m+\sigma}^{\hdot}$ en utilisant $\OmX \simeq O_{\mathbb{P}^{n+1}}(n+2)$.

La premi\`ere \'etape permet alors de conclure. 

\end{prv}

\begin{lem}\label{le33}
 Si la fl\`eche d'\'evaluation $\text{ev}_{m} : H^0(\OP(m)) \to H^0(\OS(m))$ est surjective, $H^{n+1}(K(G)^\bullet_{\sigma-m})=0$.
\end{lem}
\begin{prv}
 Au vu du lemme \ref{le31}, reprenant les notations de la premi\`ere \'etape de la preuve du lemme \ref{lemac}, on a
aussi surjectivit\'e de
$E^{n+2,0}_{n+3} \to \mathbb{H}^{n+2}(\mathcal{K}(G)_m^{\hdot})$. De ceci suit que $E_2^{1, n+1}=0$.
Or, par la seconde \'etape de la preuve du lemme \ref{lemac}, $E_2^{1, n+1}\simeq H^{n+1}(K(G)^\bullet_{\sigma-m})^{\vee}$. 
\end{prv}

Fixons $m_0$ tel que 
 \[\text{ev}_{m_0} : H^0(\OP(m_0)) \to H^0(\OS(m_0))\]
 soit surjective.

\begin{coro}
Pour tout $m \leq \sigma-m_0$, $H^{n+1}(K(G)^\bullet_m)=0$. 

\end{coro}

\begin{prv} 
Soit $m \leq \sigma-m_0$. La fl\`eche 
$\text{ev}_{\sigma-m} : H^0(\OP(\sigma-m)) \to H^0(\OS(\sigma-m))$ est  surjective. 
On peut donc conclure avec le lemme \ref{le33}.\end{prv}

Plus g\'en\'eralement, on a $h^{n+1}(K(G)^\bullet_m)=h^1(\mathcal{J}_{\Sigma}(\sigma-m))$. 

\begin{coro}\label{injectivite1}
Supposons  que $d-(n+2)>0$.  Supposons   que $A/J_{\ge d-(n+2)}$ est
 engendr\'e en degr\'e $d-(n+2)$ \footnote{Ce qui est garanti d\`es que $m_0\le d-(n+2)$.}
et que $I/J_{\ge d-(n+2)}$ est engendr\'e en degr\'e $d-(n+2)$. 
Si $m_0\le d-(n+2)$, 
l'application lin\'eaire induite par la multiplication
$$ I/J_d \to \mathrm{Hom} (A\slash J_{d-(n+2)}, I\slash J_{2d-(n+2)})
$$
est injective.
 \end{coro}
\begin{prv}
Soit $P$ un \'el\'ement du noyau. On a,  pour tout $Q'''\in (A/J)_{d-n-2}$,  $PQ'''=0\mod J$. 
De l\`a $P Q'''Q''=0\mod J$ pour tout $Q'' \in A/J_{\ge 0}$. 
Puisque   $A/J_{\ge d-(n+2)}$ est engendr\'e en degr\'e $d-(n+2)>0$, 
 il suit que $PQ'=0\mod J$ pour tout $Q'\in A/J_{\ge d-(n+2)}$. 
Comme $\sigma-2d+n+2=(n+2)(d-1)-2d=nd-n-2 \ge d-(n+2) $ on a 
$$\forall Q\in I/J_{d-n-2} \  \forall Q'\in I/J_{\sigma-2d +n+2}
 \quad <Q; PQ'>_{d-(n+2)}=0$$
o\`u $<-; - >_d$ d\'esigne l'accouplement de dualit\'e d\'efini par $d_{n+1}$ entre $I/J_d$ et $I/J_{\sigma-d}$
au lemme \ref{lemac}. 
Or on a, par fonctorialit\'e de la dualit\'e de Serre,   $$ <Q; PQ'>_{d-(n+2)}=<QQ'; P>_{\sigma-d}.$$ 
Puisque $I/J_{\ge d-(n+2)}$ est engendr\'e en degr\'e $d-(n+2)$, $P$ est 
orthogonal \`a l'espace 
$I/J_{\sigma- d}$ entier.  Donc $P=0\mod J$. 
\end{prv}

Une preuve plus courte du corollaire \ref{injectivite1} valable dans le cas  d'une hypersurface avec un seul n\oe ud 
nous a \'et\'e communiqu\'ee par  A. Otwinowska au moment o\`u nous finissions la
pr\'esente preuve. 

\section{ Th\'eor\`eme de Torelli local sur  les strates isosinguli\`eres}
\label{4}

Soit $n$  un entier pair strictement positif, et soit $X_f\subset \PP$ une hypersurface de
 degr\'e $d$ \`a singularit\'es isol\'ees simples. Comme dans le cas o\`u $X_f$ est lisse, 
une d\'eformation isosinguli\`ere de $X_f$ fournit une variation
 de structure de Hodge dont la tige en $X_f$ est $H^{n+1}(\PP-X_f)$.
Rappelons que $H_f:=H^{n+1}( \PP -X_f)$ est pur de poids $n+2$ par la discussion de la section \ref{sec2}
et que sa filtration de Hodge  v\'erifie 
\[0=F^{n+2} \subset F^{n+1} \subset \ldots \subset F^1= H^{n+1}(\PP-X_f, \C).
\]

L'IVHS correspondante  a \'et\'e  \'etudi\'e par Dimca et Saito \cite{DS} et
avec plus de d\'etails dans le cas nodal par ces m\^emes auteurs et Wotzlaw
 \cite{DSW}.
Dans cette section, nous extrayons de 
leur travail tous les renseignements dont nous aurons besoin
en ajoutant quelques petits points suppl\'ementaires. Ceci permet d'appliquer le résultat d'injectivité \ref{injectivite1} et d'aboutir à la démonstration du th\'eor\`eme \ref{torelli2}.

\subsection{Formule de Dimca-Saito-Wotzlaw pour les deux premiers termes de la filtration de Hodge}

Si $y$ est un point singulier de $X_f$, on note comme   \cite[section 1.1]{DSW}  $\wt{\alpha}_{X_f,y}$
 le plus petit z\'ero de la $b$-fonction de la singularit\'e. 
Dans notre cas, la singularit\'e est quasihomog\`ene et on note $w_1, ..., w_{n+1}$ 
les poids correspondants. Alors, il est connu que $\wt{\alpha}_{X_f,y} = \sum_i w_i$. 
Par exemple, pour une singularit\'e $A_1$, on a $\wt{\alpha}_{X_f,y} = (n+1)/2$. 
Un examen des \'equations des singularit\'es simples montre que l'on a toujours $\wt{\alpha}_{X_f,y} >1$. Plus précisemment,
 si $n\geq 4$, alors $\lfloor \wt{\alpha}_{X_f,y}\rfloor >1$, et si $n=2$, alors $\lfloor \wt{\alpha}_{X_f,y}\rfloor =1$.

Lorsque $n=2$ on d\'efinit un id\'eal homog\`ene $I'$ de l'anneau des polyn\^omes de $n+2$ variables  de la fa\c{c}on suivante. 
D\'efinissons d'abord comme \cite[(2.1.4)]{DSW} le faisceau d'id\'eaux
$\mathcal I'_{(1)} \subset \mathcal O_{\mathbb P^3}$ cosupport\'e aux points singuliers
 $y$ de $X_f$ par 
\[F_1 \mathcal O_{\mathbb P^3,y}(*X_f) = \mathcal I'_{(1)} \mathcal O_{\mathbb P^3}(2X_f),\]
o\`u le membre de gauche d\'esigne la filtration de Hodge du $D$-module $O_{\mathbb P^3}(*X_f)$
correspondant au Module de Hodge Mixte sous jacent \`a  $Rj_* \Q_{\mathbb{P}^3-X_f}$.
On pose alors $I'_k = \Gamma(\mathbb P^3,\mathcal I'_{(1)}(k))$ puis 
$I' = \oplus_{k\in\N}I'_k$. 

Par \cite[Theorem 2.2]{DSW}, on a : 
\begin{prop}\label{dsw}
 \[Gr^{n+1}_F H^{n+1}(\PP-X_f,\C)= A/J_{d-n-2},\]
\[ Gr^n_F H^{n+1}(\PP-X_f,\C)= \begin{cases} A/J_{2d-n-2} &\mbox{ si } n\ge 4\\
I'/J_{2d-n-2}&\mbox{ si }n=2\,\end{cases}\]
\end{prop}

Dans \cite[Lemma 1.5]{DSW} est \'enonc\'e que,  pour les surfaces nodales,
$\mathcal I'_{(1)}$ est l'id\'eal 
des fonctions qui s'annulent sur les points singuliers  de $X_f$ (avec structure r\'eduite).
Il est facile de g\'en\'eraliser: 
\begin{lem}
L'id\'eal $\mathcal I'_{(1)}$ co\"incide avec l'id\'eal de Tjurina:  $\mathcal I'_{(1)}=\mathcal{J}_{\Sigma}$.
\end{lem}
\begin{prv}
Il s'agit de voir que $F_1 \mathcal O_{\mathbb P^3,y}(*X_f)=(\frac{\partial h}{\partial x_1}\ldots \frac{\partial h}{\partial x_n}) h^{-2}  \mathcal O_{\mathbb P^3,y}$
o\`u  $h$ est
une \'equation locale de $X_f$ pr\`es de $y$
Or, d'apr\`es \cite[(1.3.2)]{DSW} (qui r\'ef\`ere \`a \cite{Saito}), la filtration de Hodge sur 
$\mathcal O_{\mathbb P^3,y}(*X_f)$ est dans ce cas donn\'ee par:
\[F_1 \mathcal O_{\mathbb P^3,y}(*X_f) = F_1\mathcal D_{\mathbb P^3,y}(h^{-1}\mathcal O),\]
o\`u $D_{\mathbb P^3,y}$ est filtr\'e par l'ordre de l'op\'erateur. En effet, suivant les notations de loc. cit. $k_0=0$
et $O_{\mathbb P^3,y}^{\ge 1}=O_{\mathbb P^3,y}$. 
\end{prv}

\begin{coro}
Avec les notations de la section \ref{32}, 
 $I'=I$.
\end{coro}

\subsection{Formule de Dimca-Saito pour le premier gradu\'e de la connexion de Gauss-Manin de $\PP-X_f$}

Il y a \'egalement une seconde filtration sur $H^{n+1}(\PP-X_f, \C)$ 
la filtration par l'ordre du p\^ole not\'ee $P^{\hdot}$ et l'on a $F^i\subset P^i$ \cite{DeDi}. 

 Le long de la strate isosinguli\`ere $S$ de $X_f$ dans l'espace 
projectif param\'etrisant les hypersurfaces de degr\'e $d$,
 la connexion de Gauss Manin  v\'erifie
$\nabla P^i\subset P^{i-1} \otimes \Omega^1_S$ du moins si $\dim P^i$ est localement constante pr\`es de $[X_f]$
ce qui est vrai sur un ouvert dense  $S'\subset S$
et \cite{DS} 
donne une formule pour $Gr_P\nabla_{\xi} $ si $\xi\in T_{[X_f]} S= I/(f)_d$
en termes de la partie libre du module de Brieskorn. 

Ceci est exploit\'e dans  \cite[Remarks 3.9]{DSW} dont nous tirons la proposition suivante: 

\begin{prop}\label{ds-ivhs}
Si $n\ge 4$ 
$Gr_F \nabla_{\xi}: Gr^{n+1}_F H_f\to Gr^n_F H_f$ s'identifie par l'isomorphisme de la proposition
 \ref{dsw} \`a  $-1$ fois la
multiplication $I/(f)_d\otimes A/J_{d-n-2}\to A/J_{2d-n-2}$ et si  $n=2$,  \`a
 $-1$ fois la multiplication $I/(f)_d\otimes A/J_{d-n-2}\to I/J_{2d-n-2}$.
 
\end{prop}
 
 \begin{prv}
 En faisant attention au fait que la notation $n$ 
ici correspond \`a ce qui est not\'e $n-1$ dans \cite{DSW}, cela r\'esulte
de \cite{DS} de la m\^eme fa\c{c}on que dans  \cite[Remarks 3.9]{DSW} au moins sur l'ouvert $S'$. On conclut
par passage \`a la limite. 
\end{prv}

\subsection{Preuve du th\'eor\`eme \ref{torelli2}}

Le corollaire \ref{injectivite1} signifie pr\'ecis\'ement que  
pour la famille isosinguli\`ere $(\PP\times Z(l)-\mathcal{X}(k)_{Z(l)})/Z(l)$, le noyau de 
$Gr_F \nabla_{\xi}: Gr^{n+1}_F H_f\to Gr^n_F H_f$
s'identifie via la proposition \ref{ds-ivhs} \`a $J_d \subset I_d$ qui est  l'espace 
 tangent de l'orbite de $[f]$ sous $PGL(n+2)$. Ceci \'etablit le th\'eor\`eme \ref{torelli2} 
pour $d\gg k>n+2$
gr\^ace aux r\'esultats de la section \ref{sec2}. 

Soyons maintenant plus pr\'ecis sur les conditions que  $d$ doit satisfaire pour
 la condition suffisante d'engendrement de la proposition \ref{ds-ivhs}.
Pla\c{c}ons nous en $f \in Z(k)$ d'id\'eal Jacobien $\mathcal{J}_{\Sigma}$.
Nous devons nous assurer que 
$H^0(\mathbb{P}^{n+1}, O(d-(n+2))) \to O(d-(n+2))\slash \mathcal{J}_{\Sigma}$ est surjectif
\footnote{Ce qui, puisque  $X=\mathbb{P}^{n+1}$ avec $n\ge 2$, implique que $H^0(\mathbb{P}^{n+1}, O(d')) \to O(d')\slash \mathcal{J}_{\Sigma}$ est surjectif pour 
$d'\ge d-(n+2)$ et donc que les r\'esultats de
de la section \ref{sec2} s'appliquent en posant $d'=d$.} et surtout que $I/J$ est engendr\'e sur l'anneau
de polyn\^omes $A$ en degr\'e $d-(n+2)$. Ce dernier point est garanti d\`es que 
$I$ est est engendr\'e sur l'anneau
de polyn\^omes $A$ en degr\'e $d-(n+2)$, c'est \`a dire d\`es que 
$\mathcal{J}_{\Sigma}(d-(n+2))$ est engendr\'e par ses sections globales.

Une condition suffisante est que $d-(n+2) \ge Reg(\mathcal{J}_{\Sigma})$ o\`u $Reg$ est la r\'egularit\'e 
de Castelnuovo-Mumford. 
Or, puisqu'il d\'efinit un sous-sch\'ema artinien,
 $\mathcal{J}_{\Sigma}$ est $m$-r\'egulier (pour $m\ge n+1$) si et seulement si
 $H^1(\mathbb{P}^{n+1}, \mathcal{J}_{\Sigma}(m-1)=0$,
c'est \`a dire si et seulement si $H^0(\mathbb{P}^{n+1}, O(m-1)) \to O(m-1)\slash \mathcal{J}_{\Sigma}$ est surjectif.

La condition suffisante obtenue est donc tout simplement la surjectivit\'e de:
$$
H^0(\mathbb{P}^{n+1}, O(d-(n+3))) \to O(d-(n+3))\slash \mathcal{J}_{\Sigma}.
$$

De sorte que la proposition \ref{ds-ivhs} s'applique d\`es que $d\ge n+3 + s_k(\mathbb{P}^{n+1})$, où $s_k$ est la quantité apparaissant dans l'énoncé du théorème \ref{torelli2}).

\section{Application \`a la conjecture de Shafarevich sur l'Uniformisation}

\subsection{Vari\'et\'es propres  sur $\wt{U}$.} Dans ce paragraphe on reprend les notations de \ref{11} et on ne suppose pas que $X=\PP$.

Soit $Z$ une variété connexe projective lisse, $f: Z\to \wt{U}$ un morphisme et  $z\in Z$ un point base. 
La repr\'esentation $\bar\rho: \pi_1(\wt{U}, x)\to O(H^n(X_g,\Q))$ construite \`a la proposition \ref{rhobar}
induit une repr\'esentation $f^*\bar\rho: \pi_1(Z, z)\to O(H^n(X_g,\Q))$. Notons $\wt{Z^{un}}\to Z$ le rev\^etement universel de $Z$, et $\wt Z^\rho:= \ker(f^*\bar\rho) \backslash \wt{Z^{un}}$ le rev\^etement topologique de $Z$  attach\'e \`a $f^*\bar \rho$,   

\begin{prop}
La variété $\wt Z^\rho$ est holomorphiquement convexe.
\end{prop}

\begin{prv}
Ceci r\'esulte des r\'esultats de \cite{Eys} modulo le fait que $\{ f^*\bar \rho \}$ 
est constructible absolu. 
Le cas pr\'esent est particuli\`erement simple et il est facile de d\'ecrire
la r\'eduction de Cartan-Remmert en termes de l'application des p\'eriodes. 

Notons que $f^*\bar \rho$ est
sous jacent \`a une VSH polarisable d\'efinie sur $Z$ qui est 
$f^*\V_{\wt{U}}$. Notons $\Gamma'=f^*\rho (\pi_1(Z,z))$. 
 L'application des p\'eriodes
attach\'ee \`a $f^*\V_{\wt{U}}$ se  rel\`eve \`a une application  $\Gamma'$-\'equivariante
$\mathcal P:\wt Z^\rho \to \mathcal{D}$. 
Comme $\Gamma'$ agit proprement discontin\^ument sur $\mathcal{D}$ puisque $\Gamma'$ est discret,
 il suit que $\mathcal P$
est propre. Considérons sa factorisation de Stein
\[ \wt Z^\rho\buildrel{\alpha}\over \to R \buildrel{\beta}\over\to \mathcal D, \]
où  $R$ est un espace complexe normal, $\alpha$ est propre surjective \`a fibres connexes, et $\beta$ finie. 

Comme il n'existe pas d'application holomorphe horizontale  $M \to \mathcal D$ o\`u $M$ est compacte complexe (\cite{Gri}), les fibres de $\alpha$ sont les sous espaces  analytiques connexes ferm\'es maximaux de $\wt Z^\rho$ et donc que $R$ n'a pas de sous-espace complexe analytique compact de dimension positive. 

Pour montrer que $R$ est de Stein,  on utilise la solution de Narasimhan du probl\`eme de Levi. 
On peut construire des fonctions $\mathcal C^{\infty}$ positives 
et exhaustives sur $\mathcal D$ dont le hessien complexe est d\'efini positif le long de la distribution
horizontale de $\mathcal D$ et il est ais\'e de les modifier pour construire
 une fonction  d'exhaustion  strictement plurisousharmonique sur $R$, voir \cite{Eys} pour plus de d\'etails. 

\end{prv}

\begin{coro}\label{corfinal}
 $\wt{Z^{un}}$ est  de Stein si l'application $\mathcal P$ est finie. 
\end{coro}
\begin{prv}
 En effet  $\wt Z^\rho \simeq R$ est de Stein et 
tout rev\^etement topologique d'une vari\'et\'e de Stein est Stein. 

\end{prv}
 
La discussion peut \^etre r\'esum\'ee ainsi:  si $f:Z\to \wt{U}$ est finie, 
le rev\^etement universel de $Z$ est de Stein sauf si
$Z$ contient une courbe $C$ telle que  $\V_{C}$ a monodromie finie.

 Notons $m:\wt{U} \to U$ 
resp.  $m':[\Gamma \backslash \mathcal D] \to \Gamma \backslash \mathcal D$
les applications canoniques vers les espace des modules des champs consid\'er\'es. Notons qu'il existe une application
holomorphe $p_{red}: U\to \Gamma \backslash \mathcal D$. 
Il est clair que, si $f:Z\to \wt{U}$ est finie,  $\mathcal P$ est finie si et seulement si $p_{red}$ est finie sur $Z_{r}:=m\circ f(Z)$.

\begin{prop}
 Si $f:Z\to \wt{U}$ est finie et si $Z_{r}$ est un espace projectif ou plus 
g\'en\'eralement a la propri\'et\'e que toute application holomorphe $Z_{r}\to M$ 
(avec $M$  un espace complexe) est constante ou finie,   $\wt{Z^{un}}$ est  de Stein. 
\end{prop}

\subsection{Cas où $X=\PP$.}

Le corollaire \ref{stein} r\'esulte du th\'eor\`eme \ref{torelli2} 
car celui-ci implique que l'application des p\'eriodes $p\circ f:Z\to \Gamma \backslash \mathcal{D}$ est finie sur son image. En effet, 
par le th\'eor\`eme \ref{torelli2},
une courbe $C\subset Z$ contract\'ee par $p\circ f$ est n\'ecessairement  dans une orbite
de $PGL(n+2)$. Or celles-ci sont affines car le groupe des automorphismes
birationnels d'une vari\'et\'e de type g\'en\'eral est fini et
 on d\'eduit qu'elles ne peuvent pas contenir de courbe compl\`ete.
Ceci implique que $P$ est finie  et le corollaire \ref{corfinal} permet de conclure.

{Philippe Eyssidieux}\\
{Institut Universitaire de France.  Universit\'e  Grenoble I. Institut Fourier.
100 rue des Maths, BP 74, 38402 Saint Martin d'H\`eres Cedex, France}\\
{philippe.eyssidieux@ujf-grenoble.fr}
{http://www-fourier.ujf-grenoble.fr/$\sim$eyssi/}\\

{Damien M\'egy}\\
Universit\'e de Lorraine. Institut \'Elie Cartan.
B.P. 70239, F-54506 Vandoeuvre-l\`es-Nancy Cedex, France.  \\
{Damien.Megy@univ-lorraine.fr} http://www.iecn.u-nancy.fr/$\sim$megy/
\end{document}